\documentclass[letterpaper, 10 pt, conference]{ieeeconf}
\usepackage{cite}
\usepackage{amsmath,amssymb,amsfonts}
\usepackage{algorithmic}
\usepackage{graphicx}
\IEEEoverridecommandlockouts
\allowdisplaybreaks

\renewcommand{\baselinestretch}{1}

\usepackage{color}
\usepackage{bbm}

\title{\LARGE \bf
Implied Constraint Satisfaction in Power System Optimization:\\
The Impacts of Load Variations
}

\author{Line A. Roald$^{1}$ and Daniel K. Molzahn$^{2}$% <-this % stops a space
\thanks{$^{1}$: University of Wisconsin--Madison, Madison, WI. roald@wisc.edu.}%
\thanks{$^{2}$: Georgia Institute of Technology, Atlanta, GA. molzahn@gatech.edu.}%
}

\begin{document}

\maketitle
\thispagestyle{empty}
\pagestyle{empty}

%%%%%%%%%%%%%%%%%%%%%%%%%%%%%%%%%%%%%%%%%%%%%%%%%%%%%%%%%%%%%%%%%%%%%%%%%%%%%%%%
\begin{abstract}

In many power system optimization problems, we observe that only a small fraction of the line flow constraints ever become active at the optimal solution, despite variations in the load profile and generation costs. This observation has far-reaching implications not only for power system optimization, but also for the practical long-term planning, operation, and control of the system. We formalize this empirical observation for problems involving the DC power flow equations.
We use a two-step constraint screening approach to identify constraints whose satisfaction is implied by other constraints in the problem, and can therefore be removed from the problem. 
For the first screening step, we derive analytical bounds that quickly identify redundancies in the flow limits on parallel lines. 
The second screening step uses optimization-based constraint screening, where we solve a (relaxed) optimization problem for each constraint to identify redundancies.
%As an initial step to accelerate the optimization-based constraint screening,
% that reduces the number of optimization problems that need to be solved, 
%we also perform a straightforward derivation to obtain an analytic screening method that quickly identifies redundancies in the flow limits on parallel lines. 
% This analytic method is used as a pre-processing step.
% While similar constraint screening approaches have been proposed previously, o
Different from existing methods, we specifically focus our approach on large ranges of load variation such that the results are valid for long periods of time, thus justifying the computational overhead required for the screening method.
% , which may exceed the computational effort of solving a single instance of the problem. 
% Considering implied constraint satisfaction for large variations of load also enables a range of new applications (e.g., monitoring only a subset of constraints in real-time operation and improving the tractability of stochastic optimization problems).
Numerical results for a wide variety of standard test cases show that even with load variations up to $\pm 100$\% of nominal loading, we are able to eliminate a significant fraction of the flow constraints. 
This large reduction in constraints may enable a range of possible applications.
%is suggests the suitability of constraint screening methods to address various applications that we summarize in this paper. 
As one illustrative example, we demonstrate the computational improvements for the unit commitment problem obtained as a result of the reduced number of constraints.%from the constraint screening.  
\end{abstract}

\vspace{-0.05in}

% \begin{document}
% %
% \title{Identification of Redundant Constraints\\in Power System Optimization Problems}
% %\title{Optimization-Based Constraint Screening\\ for Operational Problems in Power Systems}
% %
% %\titlend Trunning{Abbreviated paper title}
% % If the paper title is too long for the running head, you can set
% % an abbreviated paper title here
% %
% \author{Line A. Roald\inst{1} \and
% Daniel K. Molzahn\inst{2}}
% %
% %\authorrunning{F. Author et al.}
% \institute{University of Wisconsin - Madison, Madison, WI \and
% Argonne National Laboratory, Lemont, IL\\
% \email{roald@wisc.edu, dmolzahn@anl.gov}
% }
% %
% \maketitle              % typeset the header of the contribution
% %
% \begin{abstract}
% (15-250 words)

% \keywords{First keyword  \and Second keyword \and Another keyword.}
% \end{abstract}
%
%
%

\section{Introduction}
Many engineering applications give rise to optimization problems that are solved subject to feasibility constraints with time-varying problem parameters. A prominent example is electric grid optimization, where problems such as optimal power flow (OPF) and unit commitment (UC) are both important stand-alone problems in market clearing and operations and also form important building blocks in more complex problems such as long-term transmission grid planning. These problems include power flow equations that model the flow of electricity throughout the grid, incorporate technical constraints such as limits on generator outputs and transmission line capacities, and exhibit time-varying levels of electricity consumption and renewable energy generation. 
%{\color{red}We could include something about non-linearity, integer variables and how these problems are hard... Depending on how general we want to write this.}
Another important characteristic of these problems is that the number of transmission constraints is usually large. While these constraints must all be satisfied, only a limited number will be active at the optimal solution~\cite{Ng2018-bj}. 
Since the transmission constraints can be numerically challenging and represent a computational bottleneck, this observation can be exploited to devise more efficient solution algorithms, using methods such as, e.g., constraint generation~\cite{bienstock2014, roald2017corrective}.

In general, the fact that only a small number of constraints are active at the optimal solution does not imply that other constraints would not become active if the objective function changed. However, the structure of the power flow equations seems to imply that only a fraction of the constraints can ever be binding, independent of the cost function. 
Direct evidence of this is found in the literature on \emph{optimization-based bound tightening} for AC OPF~\cite{coffrin_tightening,sun2015}. Starting from an initial optimization problem with a given set of variable bounds, optimization-based bound tightening solves a sequence of optimization problems to find the upper and lower achievable values for, e.g., the voltage magnitudes and angle differences. 
The obtained values are then used to tighten the bounds on these quantities, which in turn improves the quality of certain convex relaxations of the AC power flow equations~\cite{sun2015,coffrin2015qc}. 
The effectiveness of these methods in obtaining tighter variable bounds implies that many of the bounds are implicitly satisfied through other constraints in the problem, such as the power flow equations and the generation constraints. %In this paper, we do consider the linear DC power flow approximation instead of the AC power flow equations, but the intuition about implied constraint satisfaction appears to hold true.

The results from optimization-based bound tightening can be used to remove redundant constraints from the problem before it is even passed to the solver. 
%We call this method \emph{optimization-based constraint screening}.
%{\color{red}(Is this discussed in this literature? I'm not sure about this for the same reasons we've discussed previously.)} 
A related idea was used in~\cite{ardakani2013} for identification of so-called \emph{umbrella constraints} (i.e., a subset of the constraints whose enforcement implies satisfaction of the remaining constraints). %Similar to optimization-based bound tightening, t
%The approach in~\cite{ardakani2013} also uses 
The umbrella constraints are found by solving a sequence of optimization problems to identify constraints that cannot be active, with various modifications for improving tractability in~\cite{ardakani2015} and~\cite{ardakani2018}. A related method is discussed in~\cite{madani2017}, with a focus on security-constrained unit commitment problems. 

Other related work uses analytical methods such as those presented in~\cite{chen2015,dmitry2018,molzahn-redundant_flow_limits} to provide bound tightening or constraint screening  results for AC OPF without explicitly solving optimization problems. 
%Another analytical method in~\cite{molzahn-redundant_flow_limits} identifies redundant constraints on the current flows and apparent power flows through parallel lines.
Additionally, a variety of presolvers developed for general optimization problems can be viewed as screening methods for eliminating redundant constraints~\cite{telgen1983,karwan1983,paulraj2010}. Although these analytical methods often have advantages in their computational speed, they are generally not as effective at identifying redundant constraints as optimization-based screening methods.

In this paper, we suggest to do a different kind of constraint screening. While existing papers focus on bound tightening~\cite{coffrin_tightening,sun2015} or constraint screening~\cite{ardakani2015,madani2017} for a particular instance of the power flow optimization problem where the load profile is given (or restricted to vary within a limited range), we focus on identifying constraints that are redundant across \emph{a large range of load variations}, representative of, e.g., the maximum and minimum yearly load. %The goal of this work is to determine whether it is practically interesting to perform screening in this setting, i.e., if a sufficiently large number of constraints would be found to be redundant. 
%the potential applicability of a constraint screening framework that remains valid despite large variations in the net load demands. 
Further, rather than considering a specific type of problem such as security-constrained OPF or UC, we focus more generally on optimization problems that use a DC power flow model.

Our screening approach is simple, and consists of two steps. First, we build on~\cite{molzahn-redundant_flow_limits} to provide analytic relationships for transmission constraints in the DC power flow model, which allows for the fast identification of redundant flow limits on any pair of parallel lines. %Use of this analytic screening method quickly allows us to screen out any such redundant constraints. %reduces the number of optimization problems that are solved in the optimization-based constraint screening method, which is particularly relevant in systems with many parallel lines.
Second, we use a constraint screening approach related to the optimization-based constraint screening in \cite{ardakani2015,madani2017}. However, while we still solve a (relaxed) optimization problem for each constraint, we consider a more general framework where the load parameters are also considered as optimization variables with a corresponding (large) load variation range. The screening results hence hold not only for one particular realization of the problem parameters, but for \emph{any} possible parameter realization in the considered load variation range. 

The main contribution of the paper is to empirically demonstrate that although considering larger load variations leads to an increased number of non-redundant transmission constraints, the fraction of non-redundant constraints in OPF and UC problems still remains small. This is a significant result, as it drastically increases the applicability of the proposed constraint screening method. It demonstrates that the computationally burdensome constraint screening method can be applied as an offline pre-processing step that will be valid for a prolonged period of time. 
Additionally, the long-term identification of redundant constraints may become an enabling step for a number of other applications, such as development of reduced models for long-term planning
%identify portions of a system that are subjected to high stress during certain operational regimes,
or determination of which constraints to monitor in a system with limited observability (e.g., distribution networks). 

The remainder of the paper is organized as follows. Section~\ref{sec:problem_formulation} describes the DC~OPF problem and the constraint screening method. 
Section~\ref{sec:applications} discusses several applications enabled by constraint screening with large load variations. 
Section~\ref{sec:results} numerically demonstrates the effectiveness of the constraint screening method for a variety of large systems and wide
operating ranges. 
Section~\ref{sec:conclusion} concludes the paper. 
% }

\section{Constraint Screening for~DC~Optimal~Power~Flow}
\label{sec:problem_formulation}
This paper considers a constraint screening method for power system optimization problems. These problems include the \emph{power flow equations} as a model of the network physics. The simplest such problem is the optimal power flow (OPF) problem, which seeks the minimum cost operating point that satisfies limits on generator outputs, line flows, etc. OPF formulations are important building blocks in almost all other power system optimization problems (e.g., security assessment, market clearing, long-term planning, etc.). 
%A prominent example of an extension to standard OPF problem is the transmission-constrained unit commitment (UC) problem, which includes (binary) generator on/off decisions and consideration of multiple, subsequent time periods. 

The constraint screening method in this paper is motivated by the observation that only a subset of the transmission lines ever seem to be used to their full capacity. This corresponds to a situation in which only a few of the transmission constraints are ever active. We investigate this empirically observed behavior for the DC~OPF, a linear version of the OPF problem based on the DC power flow approximation that has many important practical applications.

In particular, we propose a two-step screening method to identify redundant constraints in DC OPF problems:
\begin{itemize}
    \item An initial \emph{analytic screening} to remove redundant constraints on parallel lines. 
    % This screening step is based on topology and is evaluated using analytic relationships.
    \item An \emph{optimization-based screening} to remove redundant constraints based on techniques from optimization-based bound tightening.
\end{itemize}
We next formulate the DC OPF problem and then describe each of these steps in sequence.

%Then, as an example extension, we consider the DC UC problem.

%While there are different OPF problems in literature, we will consider the linear DC power flow approximation and the corresponding DC OPF problem for the purposes of this paper. Extensions to the more general, non-linear case involving the non-convex AC power flow constraints are left to future work. 
%however, our preliminary results indicate that the structure of the power system problems that make the method effective will
%type of optimization problems 
%To set the stage for the discussion about the constraint screening approach, we describe standard formulations of the DC OPF and the DC UC below. 

\subsection{DC Optimal Power Flow}
\label{subsec:dcopf}
The transmission network is modelled as a graph $(\mathcal{V,\mathcal{L}})$, where $\mathcal{V}$ denotes the nodes and $\mathcal{L}$ denotes transmission lines. 
%Each line $(i,j) \in \mathcal L$ is given an (arbitrary) orientation, with the convention that the power flowing from $i$ to $j$ is positive, while the power flowing from $j$ to $i$ is negative. 
The numbers of nodes and lines are $|\mathcal{V}| = v$ and $|\mathcal{L}|=\ell$, respectively. The set of dispatchable generators is denoted by $\mathcal{G} \subseteq \mathcal V$, and the total number of generators is $|\mathcal{G}|=g$. 
%The vector $b\in \R^\ell$ denotes the transmission line susceptances, and 
The vector $d$ encodes the demands at each bus.
The decision variables are the active power generated by each generator $p_i$ for all $i \in \mathcal G$ and the voltage angles $\theta_i$ for all $i \in \mathcal V$.
The maximum power flow across each transmission line $ij\in\mathcal{L}$ is given by $f_{ij}^{\text{max}}$, and the generator limits are given by $p_i^\text{min}$ and $p_i^\text{max}$.
For ease of notation and without loss of generality, we will assume that there is one generator per node, such that $\mathcal{G}=\mathcal V$ and $g=n$. This is easily extended to the more general case with zero or multiple generators per node by including a matrix which maps each generator to its respective node. In addition, to keep notation clear, we assume that $p_i^\text{min},p_i^\text{max} \geq 0$, although extensions are straightforward.

The DC OPF problem finds the optimal generation dispatch $p^*$ which minimizes operational cost:

\begin{subequations} \label{eq:opf}
\begin{align}
    \underset{p, \theta}{\text{min }}\ &\sum_{i\in\mathcal{G}} c_{0,i} + c_{1,i} p_i + c_{2,i} p_i^2&
    \label{eq:objective}\\
    \text{s.t. }\ &p_i - d_i = \sum_{j:(i,j)\in\mathcal{L}} b_{ij}(\theta_i - \theta_{j}),
    &&\forall ~i\in\mathcal{V}, \label{eq:balance}\\
	&p_i^{\text{min}}\leq p_i\leq p_i^{\text{max}}, 
	&&\forall ~~i\in\mathcal{G},\label{eq:generation}\\
	-&f_{ij}^{\text{max}}\leq b_{ij}(\theta_i - \theta_{j})\leq f_{ij}^{\text{max}},
	&&\forall ~ij\in\mathcal{L}.
	\label{eq:transmission}
\end{align}
\end{subequations}
The objective~\eqref{eq:objective} minimizes the generation cost, modelled as a function with constant, linear, and quadratic cost coefficients $c_{0,i}$, $c_{1,i}$, and $c_{2,i}$, $\forall i \in\mathcal{G}$. Equation~\eqref{eq:balance} enforces power balance at every node. 
The power flow from node~$i$ to each adjacent node~$j$ is modelled as a function of the angle difference $\theta_i - \theta_j$ and the susceptance $b_{ij}$.
Limits on generator outputs and line flows are enforced by~\eqref{eq:generation} and~\eqref{eq:transmission}. 

Constraint screening identifies redundant constraints among the line flow limits~\eqref{eq:transmission}. The set of redundant constraints is denoted as $\mathcal{R}\subseteq\mathcal{M}$, where $\mathcal{M}$ is the set of upper and lower bounds on the transmission lines \eqref{eq:transmission}, and $|\mathcal{M}|=2|\mathcal{L}|=2\ell$. There is no need to perform screening for the power balance constraints \eqref{eq:balance} as equality constraints will never be redundant. Further, note that we do not screen for the generator inequality constraints for two reasons: First,  we observed in experiments not reported here that the number of generator constraints that can be screened out is very small. This is as expected, as generators are typically located at nodes which can accommodate the full ranges of the generators' power outputs. Second, the generation limits are bounds on a single variable, which are generally easy to handle numerically. In contrast, the transmission line constraints involve a combination of several variables
% \footnote{This is particularly easy to observe in an equivalent and commonly used OPF formulation based on the so-called Power Transfer Distribution Factors (PTDFs) \cite{}, where the power flows are expressed as direct functions of the generator power injections}, 
and can represent a computational bottleneck.

\subsection{Constraint Screening for DC Optimal Power Flow}
\label{sec:all_screening}

\subsection*{Step 1 -- Analytic Constraint Screening: \\ Implied Constraint Satisfaction on Parallel Lines}
\label{sec:par}
The first step in our constraint screening is an analytic analysis that removes redundant constraints on parallel lines. 
%This method is based on the topology of the system and is independent of the operating point. 
While the derivation of this step is straightforward, experiments show its ability to quickly identify a non-negligible number of redundant constraints for certain test cases.
% {\color{red}We might want to clarify quite clearly what we are adopting from literature and what is new... I don't know in which category this belongs... It is probably fundamental enough to be textbook material, but we don't know where to find it written out??}

Parallel lines are lines which share the same terminal nodes. While parallel lines may have different impedances and flow limits, they have the same set of terminal voltages. Thus, the flows on parallel lines are not independent and enforcing the flow limit for one line may imply satisfaction of the limits for other parallel lines. Prior work~\cite{molzahn-redundant_flow_limits} developed a condition for redundancy of the flow limits for parallel lines using an AC power flow model. This section derives an analogous condition using a DC power flow model.

Consider the set of parallel lines connected between buses~$i$ and~$j$, with superscripts $(\,\cdot\,)^{(k)}$ and $(\,\cdot\,)^{(l)}$ denoting quantities associated with an arbitrary pair of these parallel lines. 
Using the DC power flow model, the line flows are
\begin{subequations}
\begin{align}
 f_{ij}^{(k)} = b_{ij}^{(k)} \left( \theta_i - \theta_j\right), \quad
 f_{ij}^{(l)} = b_{ij}^{(l)} \left( \theta_i - \theta_j\right). \nonumber
\end{align}
\end{subequations}
Observe that all the susceptances $b_{ij}^{(k)}$ and $b_{ij}^{(l)}$ are positive, and thus the power flows on two parallel lines are always in the same direction. In addition, the flow limits $f_{ij}^{(k),\max}$, $f_{ij}^{(l),\max}$ are also positive, and the upper and lower limits on the flows in~\eqref{eq:transmission} are symmetric, i.e., 
$f_{ij}^{(k),\min} = -f_{ij}^{(l),\max}$. Hence, we only need to consider the absolute value of the angle difference, $\left| \theta_i - \theta_j\right|$. %In addition, the flow limits $f_{ij}^{(k),\max}$, $f_{ij}^{(l),\max}$ are also positive, and the upper and lower limits on the flows in~\eqref{eq:transmission} are symmetric, i.e., the lower limits are the negative of the upper limits.

The flow constraint on the $k^{th}$ parallel line is redundant if the power flow on the $l^{th}$ parallel line always reaches its limit $f_{ij}^{(l),\max}$ prior to the power flow on the $k^{th}$ parallel line reaching its limit $f_{ij}^{(k),\max}$. 
To analyze whether this condition holds, we normalize the line constraints by their limits, i.e.,
\begin{equation}
    \frac{b_{ij}^{(k)} \left| \theta_i - \theta_j\right|}{f_{ij}^{(k),\max}}\leq 1, \qquad \frac{b_{ij}^{(l)} \left| \theta_i - \theta_j\right|}{f_{ij}^{(l),\max}} \leq 1. \nonumber
    %, & \forall \theta_i,\,\theta_j
\end{equation}
To show that one constraint reaches its limit before the other, it now suffices to show that
\begin{align}
    \label{par_condition_orig}
    & \frac{b_{ij}^{(k)} \left| \theta_i - \theta_j\right|}{f_{ij}^{(k),\max}} < \frac{b_{ij}^{(l)} \left| \theta_i - \theta_j\right|}{f_{ij}^{(l),\max}}
\end{align}
for any operating condition, i.e., independent of the specific values for $\theta_i,\,\theta_j$.
%We next simplify~\eqref{par_condition_orig} by
However, since the two lines share the same voltage angles $\theta_i,~\theta_j$ at the terminals, we can simplify \eqref{par_condition_orig}:
\begin{align}\label{par_redundant_1}
    \frac{b_{ij}^{(k)}}{f_{ij}^{(k),\max}} < \frac{b_{ij}^{(l)}}{f_{ij}^{(l),\max}}.
\end{align}
%
% Similarly, the flow limit on the $l^{th}$ parallel line between buses~$i$ and~$j$ is redundant if the following condition holds:
% %
% \begin{align}\label{par_redundant_2}
%     \frac{b_{ij}^{(k)}}{f_{ij}^{(k),\max}} > \frac{b_{ij}^{(l)}}{f_{ij}^{(l),\max}}.
% \end{align}
% %
Satisfaction of \eqref{par_redundant_1} guarantees redundancy of the flow limit on the $k^{th}$ parallel line between buses~$i$ and~$j$
Thus, to identify redundant constraints on parallel lines, we evaluate  $b_{ij}^{(l)}/f_{ij}^{(l),\max}$ for all the lines and check~\eqref{par_redundant_1}  %and~\eqref{par_redundant_2} 
for all pairs of parallel lines. We only keep the constraint with the larger fraction (i.e., line $(l)$ in \eqref{par_redundant_1}) which will reach its limit first, and add all other constraints to the set of redundant constraints $\mathcal{R}$.
%provides a fast method for screening out unnecessary flow limits in the DC power flow model. 
Note that equality of the left and right sides of~\eqref{par_redundant_1} indicates that the constraints are the equivalent and thus either can be arbitrarily selected as redundant. Observe that this screening method only depends on the system topology and parameters, with no dependence on the operating point. The constraints identified as redundant will hence remain redundant regardless of the load and generation profiles.
%independent and thus only needs to be computed once.

\subsection*{Step 2 -- Optimization-based Constraint Screening: \\ Implied Constraint Satisfaction over Varying Loads}
\label{sec:screening}

The second step applies an optimization-based constraint screening method inspired by optimization-based bound tightening~\cite{coffrin_tightening,sun2015} that is closely related to previous methods for identifying so-called ``umbrella'' constraints~\cite{ardakani2013,ardakani2015,madani2017,ardakani2018}. %These methods solve optimization problems in order to identify bounds on certain quantities that are implied by other constraints in the problem.
% , and have been shown to be very effective to generate better variable bounds in AC OPF \cite{}. Here, we apply a similar procedure to the DC OPF.
% The set-up is similar to \cite{}, but is less focused on unit commitment and considers much larger ranges of loads. 
However, while the our second constraint screening step is similar to previously proposed methods, the goal of our screening is to consider large load variations. 
% {\color{red}The ability to obtain effective screening results (see Section~\ref{sec:results}) despite considering large ranges of variation is crucial for success in several of our proposed applications discussed in Section~\ref{sec:applications}.}

\subsubsection{Per-Constraint Optimization Problem}
The key idea of the optimization-based screening method is to solve a modified version of the original optimization problem for each transmission line constraint \eqref{eq:transmission} to identify whether or not the constraint can ever be active. 
%The parameters $\xi$ are treated as optimization variables that fall within a specified uncertainty set. 
%The optimization problems solved for each constraints are modified versions of the original optimization problem. 
%These modified problems are \emph{relaxations} of the original optimization problem, and 
In those modified problems, the original objective function is replaced by an objective which maximizes or minimizes the value of the power flow on the transmission line under consideration. 
The set of decision variables still includes the original optimization variables $p$ and $\theta$ (with their original bounds). Since we are interested in certifying constraint redundancy not only for one load profile, but for a larger range of load variations, we also include the load demands $d$ as optimization variables in the modified problem. The loads can take any value within a predetermined polyhedral uncertainty set $\mathcal{D}$. 

This gives rise to the following optimization problems to obtain the maximum (and minimum) achievable power flow $f_{mn}$ for each line $mn \in\mathcal{L}$,
\begin{subequations} \label{eq:opf_mod}
\begin{align}
    &\underset{p, \theta, d}{\text{max}} &&\!\!\!\!\!\!\!\!\!/\,\underset{p, \theta, d}{\text{min}}~\ f_{mn} &
    \label{eq:objective_mod}\\
    &\text{s.t. }\ &&f_{mn}=b_{mn}(\theta_m - \theta_{n}),\\
    %& \text{\eqref{eq:balance},\eqref{eq:generation}, \eqref{eq:transmission} }
    &&&p_i - d_i = \sum_{j:(i,j)\in\mathcal{L}} b_{ij}(\theta_i - \theta_{j}),
    &&\forall ~i\in\mathcal{V}, \label{eq:balance_mod}\\
	&&&p_i^{\text{min}}\leq p_i\leq p_i^{\text{max}} ,
	&&\forall ~~i\in\mathcal{G},\label{eq:generation_mod}\\
	&&-&f_{ij}^{\text{max}}\leq b_{ij}(\theta_i - \theta_{j})\leq f_{ij}^{\text{max}},
	&&\forall ~ij\in\mathcal{L},
	\label{eq:transmission_mod}\\
	&&&d\in\mathcal{D}. && \label{eq:uncertainty_mod}
\end{align}
\end{subequations}
%minimization or maximization problem for each constraint, depending on whether the constraint is an lower or upper bound, respectively. 

\subsubsection{Certificate of Implied Constraint Satisfaction}
If the maximum (or minimum) power flow given by the objective function value $f_{mn}^*$ does not achieve the constraint bound, i.e., $f_{mn}^*< f_{mn}^{\text{max}}$ for the maximization problem (or $f_{mn}^*> -f_{mn}^{\text{max}}$ for the minimization problem), we have a certificate that the upper (or lower) bound on the transmission line $mn$ in \eqref{eq:transmission} can never be violated for any load variation described by $d\in\mathcal{D}$. In this case, we can guarantee that the constraint will be satisfied even if it is not explicitly considered in the model and we add the constraint to the set of redundant constraints, which we denote by $\mathcal{R}$.

Solving~\eqref{eq:opf_mod} for each constraint is computationally expensive, particularly in systems with a large number of transmission lines $\ell$, as each problem has similar complexity to the original OPF problem. Therefore, we are interested in obtaining screening results that consider a large set of possible load profiles $\mathcal{D}$, such that we do not need to rerun the screening frequently. 
However, while we would like the results to be valid for prolonged periods of time, considering a larger set $\mathcal{D}$ increases the feasible set of~\eqref{eq:opf_mod}, leading to a larger range of achievable power flows $f_{mn}^*$. This again implies that fewer constraints will be deemed redundant, hence making the screening procedure less effective.

\subsubsection{Reduced DC OPF After Screening}
Given the set of redundant flow constraints $\mathcal{R}$, we solve the following reduced DC OPF problem:
\begin{subequations} \label{eq:opf_red}
\begin{align}
    \underset{p, \theta}{\text{min }}\ &\sum_{i\in\mathcal{G}} c_{0,i} + c_{1,i} p_i + c_{2,i} p_i^2&
    \label{eq:objective_red}\\
    \text{s.t. }\ &\eqref{eq:balance},~\eqref{eq:generation},
	&&\label{eq:summary}\\
	-&f_{ij}^{\text{max}}\leq b_{ij}(\theta_i - \theta_{j})\leq f_{ij}^{\text{max}},
	&&\forall ~ij\in\mathcal{L}\backslash\mathcal{R},
	\label{eq:transmission_red}
\end{align}
\end{subequations}
where we consider a smaller number of transmission constraints $\mathcal{L}\backslash\mathcal{R}$.

\subsubsection{DC OPF Reformulation for Problems with Few Line Constraints}
The DC~OPF problem~\eqref{eq:opf} can be equivalently formulated using so-called Power Transfer Distribution Factors (PTDFs)~\cite{wood2013}. This formulation eliminates the $\theta$ variables by substitution through the nodal power balance equations~\eqref{eq:balance}, and directly expresses the line flow constraints as a function of the power injections. Given the PTDF-formulation of the DC OPF, the reduced problem \eqref{eq:opf_red} can be equivalently expressed as 
\begin{subequations} \label{eq:ptdf_opf}
\begin{align}
    \underset{p, \theta}{\text{min }}\ &\sum_{i\in\mathcal{G}} c_{0,i} + c_{1,i} p_i + c_{2,i} p_i^2&
    \label{eq:ptdf_objective}\\
    \text{s.t. }\ &\sum_{i\in\mathcal{V}}p_i - d_i =0, \label{eq:ptdf_balance}\\
	&p_i^{\text{min}}\leq p_i\leq p_i^{\text{max}}, 
	&&\forall ~~i\in\mathcal{G},\label{eq:ptdf_generation}\\
	-&f_{ij}^{\text{max}}\leq \mathbf{M}_{(ij,\cdot)}(p-d) \leq f_{ij}^{\text{max}},
	&&\forall ~ij\in\mathcal{L}\backslash\mathcal{R}.
	\label{eq:ptdf_transmission}
\end{align}
\end{subequations}
Here, the matrix $\mathbf{M}$ is referred to as the PTDF matrix and is defined in~\cite{wood2013}. The notation $\mathbf{M}_{(ij,\cdot)}$ indicates the row of $\mathbf{M}$ corresponding to the line $ij \in \mathcal{L}$. 

The main difference between the $\theta$-formulation \eqref{eq:opf_red} and the PTDF-formulation \eqref{eq:ptdf_opf} can be summarized based on the number of variables and the level of sparsity. The $\theta$-formulation has a larger number of variables, but a sparse set of constraints. The PTDF-formulation, on the other hand, has a smaller number of variables, but the constraint matrix $\mathbf{M}$ is dense. Which formulation is better hence usually depends on the ability of a given solver to handle a larger number of variables and exploit sparsity. However, when the set of non-redundant constraints $\mathcal{L}\backslash\mathcal{R}$ is small, the PTDF formulation has the advantage of considering a much smaller number of variables and only a few dense constraints.

\subsection{Constraint Screening for More General Power System Optimization Problems}
\label{subsec:more_general_screening}
The constraint screening approach described above can be extended to optimization problems that involve the DC power flow constraints as part of a more complex problem formulation including, e.g., integer variables, consideration of multiple time periods, and non-linear constraints. The idea is straightforward: instead of doing constraint screening directly on the more complex problem (which would involve solving the complex problem a large number of times), we \emph{relax} the complex optimization problem to a simpler problem which resembles the linear DC OPF formulation. 

\subsection*{Illustrative Example: DC Unit Commitment} 

To illustrate the idea of applying the screening to more complex optimization problems, we consider the so-called DC Unit Commitment (DC UC) problem. The DC UC is an extension of \eqref{eq:opf} that accounts for decisions related to the start-up and shut-down of generators. This is particularly important for scheduling in problems with non-zero no-load cost $c_0$ and non-zero lower generation bounds $p_i^{\text{min}}$. In such cases, consideration of generator shut-down may lead to more economical solutions and might be necessary to obtain feasibility, e.g., if the demanded power falls below the minimum generator bounds. Since the on/off decisions are naturally binary, the DC UC problem is a mixed-integer linear program (MILP). It hence belongs to a significantly harder class of optimization problems than the DC OPF, which is a linear program (LP). This difficulty is generally reflected in longer solution times. The computational burden of the UC problem can be significant, which motivates the application of constraint screening methods. However, to make the screening problems efficiently solvable, we first relax the UC problem to bring it to a simpler form, which is easier to handle computationally. The analytical first step of our screening method can be applied directly without modification to remove redundant flow constraints on parallel lines. Modifications to the optimization-based second step are described in the remainder of this section.

We note that previous work in~\cite{ardakani2015} and~\cite{madani2017} also considers constraint screening for DC UC problems. The focus of our constraint screening work differs from this prior work in that we specifically consider the validity of wide ranges of load variation with the goal of facilitating extensions to various applications discussed in Section~\ref{sec:applications}.

\subsubsection{Relaxation of Multi-Period Constraints}
The generator start-up and shut-down decisions considered in the DC UC problem generally require the consideration of multiple time periods, as any given generator has limitations on its minimum up-time and minimum down-time after the generator is started up and shut down, respectively. Other time-coupling characteristics include start-up costs that depend on how long the generator has been turned off and ramping constraints that restrict the ability of a generator to adjust its set-points between periods. See~\cite{knueven2018} for an extensive model of a typical DC~UC problem.

As a first step to simplify the UC model, we relax any multi-period constraints by simply removing them from the problem. This leads to a set of decoupled unit commitment problems which only involve a single time period OPF problem and additional binary decision variables $z_i$ for all $i \in \mathcal G$ that model whether each generator is on ($z_i = 1$) or off ($z_i=0$). 
%we consider a stylized version of the UC problem which only involves a single time period and binary decision variables $z_i$ for all $i \in \mathcal G$ that model whether each generator is on ($z_i = 1$) or off ($z_i=0$). 
With this, the single-period UC problem is
\begin{subequations} \label{eq:uc}
\begin{align}
    \underset{p, \theta}{\text{min }}\ &\sum_{i\in\mathcal{G}} c_{0,i}z_i + c_{1,i} p_i + c_{2,i} p_i^2&
    \label{eq:uc_objective}\\
    \text{s.t. }\ &p_i - d_i = \sum_{j:(i,j)\in\mathcal{L}} b_{ij}(\theta_i - \theta_{j}),
    &&\forall ~i\in\mathcal{V}, \label{eq:uc_balance}\\
	&p_i^{\text{min}}z_i\leq p_i\leq p_i^{\text{max}}z_i, 
	&&\forall ~~i\in\mathcal{G},\label{eq:uc_generation}\\
	-&f_{ij}^{\text{max}}\leq b_{ij}(\theta_i - \theta_{j})\leq f_{ij}^{\text{max}},
	&&\forall ~ij\in\mathcal{L},
	\label{eq:uc_transmission}\\
	&z_i\in\{0,1\}, &&\forall~~i\in\mathcal{G}. \label{eq:integral}
\end{align}
\end{subequations}
This problem represents a relaxation of the more comprehensive UC problem, as the additional constraints involving the coupling of subsequent time-steps have been removed. However, it still includes another major complication associated with the UC problem, namely the binary variables~$z_i$.

\subsubsection{Relaxation of Integer Variables}
\label{subsubsec:relax_int}
%Applying the optimization-based constraint screening directly to the UC problem would still cause significant computational overhead, as problems involving integer variables are generally harder to solve. Therefore, t
To obtain a more tractable problem without integer variables, we relax problem \eqref{eq:uc} to a linear program. In particular, we modify \eqref{eq:uc_generation} by setting $z_i=0$ in the lower bound and $z_i=1$ in the upper bound, which corresponds to relaxing the lower generation bound by setting it to zero, i.e., $p_i^{\text{min}}=0$.
%we first relax the lower limits on the generators' power outputs to include zero. 
Since our constraint screening approach does not consider the objective function, but only investigates constraint redundancy based on the feasible region of the problem, we do not consider how the change in $z_i$ affects the objective function \eqref{eq:uc_objective}.

Hence, the relaxed problem for the UC screening is
\begin{subequations} \label{eq:opf_mod_uc}
\begin{align}
    &\underset{p, \theta, d}{\text{max}} &&\!\!\!\!\!\!\!\!\!/\,\underset{p, \theta, d}{\text{min}}~\ f_{mn}^{\text{(R)}} &
    \label{eq:objective_mod_uc}\\
    &\text{s.t. }\ &&f_{mn}^{\text{(R)}}=b_{mn}(\theta_m - \theta_{n}),\\
    %& \text{\eqref{eq:balance},\eqref{eq:generation}, \eqref{eq:transmission} }
    &&&p_i - d_i = \sum_{j:(i,j)\in\mathcal{L}} b_{ij}(\theta_i - \theta_{j}),
    &&\forall ~i\in\mathcal{V}, \label{eq:balance_mod_uc}\\
	&&&0\leq p_i\leq p_i^{\text{max}} ,
	&&\forall ~~i\in\mathcal{G},\label{eq:generation_mod_uc}\\
	&&-&f_{ij}^{\text{max}}\leq b_{ij}(\theta_i - \theta_{j})\leq f_{ij}^{\text{max}},
	&&\forall ~ij\in\mathcal{L},
	\label{eq:transmission_mod_uc}\\
	&&&d\in\mathcal{D}. && \label{eq:uncertainty_mod_uc}
\end{align}
\end{subequations}
which is similar to \eqref{eq:opf_mod} except for the relaxed lower bound $p_i^{\text{min}}=0$ in \eqref{eq:generation_mod_uc}.

\subsubsection{Certification of Implied Constraint Satisfaction}
Similar to the constraint screening method in \eqref{eq:opf_mod}, we consider cases where the maximum (or minimum) power flow given by the objective function value $f_{mn}^{\text{(R)}*}$ does not achieve the constraint bound, i.e., $f_{mn}^{\text{(R)}*}< f_{mn}^{\text{max}}$ for the maximization problem (or $f_{mn}^{\text{(R)}*}> -f_{mn}^{\text{max}}$ for the minimization problem).
In this case, we know that the power flow on transmission line $mn$ can never achieve its limit as long as the load stays within the predefined set $\mathcal{D}$ and the other constraints in \eqref{eq:opf_mod_uc} are satisfied. This provides a certificate for redundancy of the constraint on $f_{mn}$ in \eqref{eq:opf_mod_uc}. However, since \eqref{eq:opf_mod_uc} is a relaxed (less restrictive) version of the UC problem \eqref{eq:opf_mod_uc}, the objective value $f_{mn}^{*\text{(R)}}$ of the relaxed problem \eqref{eq:opf_mod_uc} is an upper bound on the objective value $f_{mn}^{*\text{(UC)}}$ of the original problem \eqref{eq:uc} in the case of maximization, or a lower bound in the case of minimization. Hence, we have that 
\begin{align}
 &f_{mn}^{*\text{(R)}}<f_{mn}^{\text{max}} \implies     f_{mn}^{*\text{(UC)}}\leq f_{mn}^{*\text{(R)}}<f_{mn}^{\text{max}}~, \\
  &f_{mn}^{*\text{(R)}}>f_{mn}^{\text{min}} \implies     f_{mn}^{*\text{(UC)}}\geq f_{mn}^{*\text{(R)}}>f_{mn}^{\text{min}}~,
\end{align}
for the upper and lower bounds, respectively.
Therefore, if for any transmission line constraint is deemed redundant by solving \eqref{eq:opf_mod_uc}, we can safely add the corresponding transmission constraint in \eqref{eq:uc} to the set $\mathcal{R}_{\text{UC}}$ of redundant constraints.

Since we are working with a relaxation, the set of certified redundant constraints $\mathcal{R}_{\text{UC}}$ is a subset of the true set of redundant constraints. This is easily seen by considering the case where
\begin{equation}
    f_{mn}^{*\text{(UC)}}< f_{mn}^{*\text{(R)}}=f_{mn}^{\text{max}}~.
\end{equation}
This constraint would not be certified as redundant by our screening method based on the relaxed problem, but would be redundant for the original UC problem since $f_{mn}^{*\text{(UC)}}<f_{mn}^{\text{max}}$.
However, while relaxing the UC problem may reduce the number of identified redundant constraint, formulating our screening problems as LPs rather than MILPs and reducing their size provides a substantial computational benefit. As will be demonstrated in the numerical test summarized in Section~\ref{sec:results}, the screening method is capable of identifying many redundant constraints despite this relaxation.

\subsubsection{Reduced UC Problem after Screening}
Similar to the DC OPF problem, we can express the reduced UC problem in terms of the PTDF matrix:
\begin{subequations} \label{eq:ptdf_uc}
\begin{align}
    \underset{p, \theta}{\text{min }}\ &\sum_{i\in\mathcal{G}} c_{0,i}z_i + c_{1,i} p_i + c_{2,i} p_i^2&
    \label{eq:ptdf_objective_uc}\\
    \text{s.t. }\ &\sum_{i\in\mathcal{V}}p_i - d_i =0, \label{eq:ptdf_balance_uc}\\
	&p_i^{\text{min}}z_i\leq p_i\leq p_i^{\text{max}}z_i, 
	&&\forall ~~i\in\mathcal{G},\label{eq:ptdf_generation_uc}\\
	-&f_{ij}^{\text{max}}\leq \mathbf{M}_{(ij,\cdot)}(p-d) \leq f_{ij}^{\text{max}},
	&&\forall ~ij\in\mathcal{L}\backslash\mathcal{R}.
	\label{eq:ptdf_transmission_uc}
\end{align}
\end{subequations}
For simplicity, we focus on the single-period UC problem.

\section{Applications Enabled by Implied Constraint Satisfaction across Large Ranges of Load}
\label{sec:applications}
The constraint screening literature~\cite{ardakani2013,ardakani2015,ardakani2018,madani2017,molzahn-redundant_flow_limits} typically focuses on the advantages of these methods with respect to the computational speed of a particular problem instance (i.e., in the context of a problem which will be solved now or in the near future).
We specifically consider much larger load variations that are representative of, e.g., all possible load profiles over a substantially longer period. 
%While computational advantages of a particular problem can indeed be an important consideration, our setup specifically considers much larger load variations that are representative of, e.g., all possible load profiles over a substantially longer period. The numerical results in Section~\ref{sec:results} show that the screening is still effective in this setting. %which enables a number of new applications. %However, there are also a variety of other applications of constraint screening that have received less attention in prior literature but are neverthelesss important for many practical purposes. 
The ability to identify redundant constraints despite large ranges of load variation enables a variety of other applications. While we have already mentioned a few applications earlier in the paper, we provide an overview here. 
%means that constraint screening methods hold significant promise for a variety of other applications that have received less attention in prior literature but are nevertheless important for many practical purposes. This section reviews a few such applications. 
%and to lay out an agenda for future research.

% This section first briefly discusses these applications in order to lay out an agenda for future research. This section then formulates the DC UC problem that serves as an illustrative test case for the computational advantages available from constraint screening methods despite large ranges of load variation.

%Simplifications:
 \emph{1) Removing constraints in operational problems:} Although the screening process is computationally expensive, the offline screening may reduce the time of solving operational problems in real-time, e.g., in market clearing based on DC OPF and DC UC with contingency constraints~\cite{ardakani2013,ardakani2015,ardakani2018,madani2017}. %Thi a sufficiently large load range. %This can also be useful for, e.g., long-term planning studies, where the operators solve a sequence of problems representative of several years of operations.

 \emph{2) Removing constraints from chance-constrained and robust DC OPF problems:} 
Many previously proposed DC OPF formulations consider joint chance constraints \cite{vrakopoulou2013} or robust constraints \cite{warrington2013} to guarantee security against uncertain load variations. Knowing that some limits will be violated before others for the considered range of uncertain loads allows for a reduction in the number of constraints in such problems.

\emph{3) Removing constraints in long-term planning problems:} Long-term planning problems can be challenging to solve due to the consideration of a very large number of time periods with significant variations in load and, more recently, renewable energy generation~\cite{lumbreras2016}. By using the constraint screening as a pre-processing step, we may be able to significantly improve the tractability of these problems. 

\emph{4) Constructing reduced models:} If we can conclusively determine that certain constraints are not required in the OPF problem, we may be able to not only remove the constraints from the OPF problem itself, but also simplify the construction of reduced models which still capture all important constraints, such as the method in \cite{jang2013}.

\emph{5) Guarantee safe operations with limited measurements and sensing:}
Operators of some power systems, such as distribution networks, have only limited real-time observability of the system state. In these systems, 
constraint screening can be used to distinguish between constraints which are not important to observe in real time (redundant constraints), and constraints which should be monitored as they will be violated before other constraints (non-redundant constraints). This type of application, which is further discussed in~\cite{hicss2019}, may enable online control applications such as~\cite{gan2016,colombino2018} by reducing the measurement requirements.

\section{Numerical Results}
\label{sec:results}

This section numerically demonstrates the capabilities of the constraint screening method described in Section~\mbox{\ref{sec:all_screening}} for a diverse set of large test cases. The constraint screening method is implemented in Julia using \mbox{PowerModelsAnnex.jl}~\cite{Coffrin2018} and Ipopt~\cite{wachter2006implementation}. The single-period UC problem~\eqref{eq:ptdf_uc} is implemented using YALMIP~\cite{yalmip} and solved with Mosek~\mbox{v.9.0.79} on a laptop computer with a quad-core 2.70~GHz processor and 16~GB of RAM. The test cases are taken from PGLib~v.17.08~\cite{pglib}. For the sake of brevity, we present detailed results for a selected subset of the test cases whose characteristics are summarized in Table~\ref{t:case_summary}. Note that many of the generators in the PGLib test cases have lower generation limits set to zero. To obtain more challenging UC problems, we modified the generators' lower limits to be the maximum of the value specified in the dataset and $10$\% of their upper generation limit. We consider a range of load demands from $0$\% to $\pm 100$\% variation around the nominal values specified in the datasets. This gives us the following description of the set of considered load ranges~$\mathcal{D}$,
\begin{equation}
    \mathcal{D} = \{(1-v)\cdot d_i^{nom} \leq d \leq (1+v)d_i^{nom}, \quad \forall i\in\mathcal{V}\}, \label{eq:uncset}
\end{equation}
where we consider $v=\{0, 0.25, 0.5, 0.75, 1\}$.

\subsection{Constraint Screening for Varying Ranges of Load}
% 1) Constraint screening results for pglib-opf cases for varying ranges of parameter variations (and relaxed generation lower bounds)

This section presents numerical results regarding the effectiveness of the constraint screening for large ranges of load variation. We first applied the analytic step described in Section~\ref{sec:par} to identify redundant flow constraints on parallel lines. This step is fast, with computation times taking less than $0.5$ seconds for all considered systems with less than 4000 buses and $4.6$ seconds for the largest considered test case case9241\_pegase. We then applied the optimization-based constraint screening method described in Section~\mbox{\ref{sec:screening}} to the remaining line flow constraints. As described in Section~\ref{subsubsec:relax_int}, the lower generation limits are relaxed to zero in order to accommodate the ability to shut off generators in the UC problem considered later in this section. The resulting linear programs are solved quickly, using on average less than 0.6s per optimization problem even for the largest test case case9241\_pegase. Due to the large number of constraints, the time required for constraint screening is still significant if the problems \eqref{eq:opf_mod} are solved sequentially. However, since each of the problems \eqref{eq:opf_mod} are completely decoupled, they can easily be computed in parallel.

%{\color{red}(DKM: We should say something about the constraint screening time required for this step here.)}

\begin{figure*}[ht]
	\centering
	\includegraphics[width=1.05\textwidth,trim={3.25cm 0cm 1cm 0},clip]{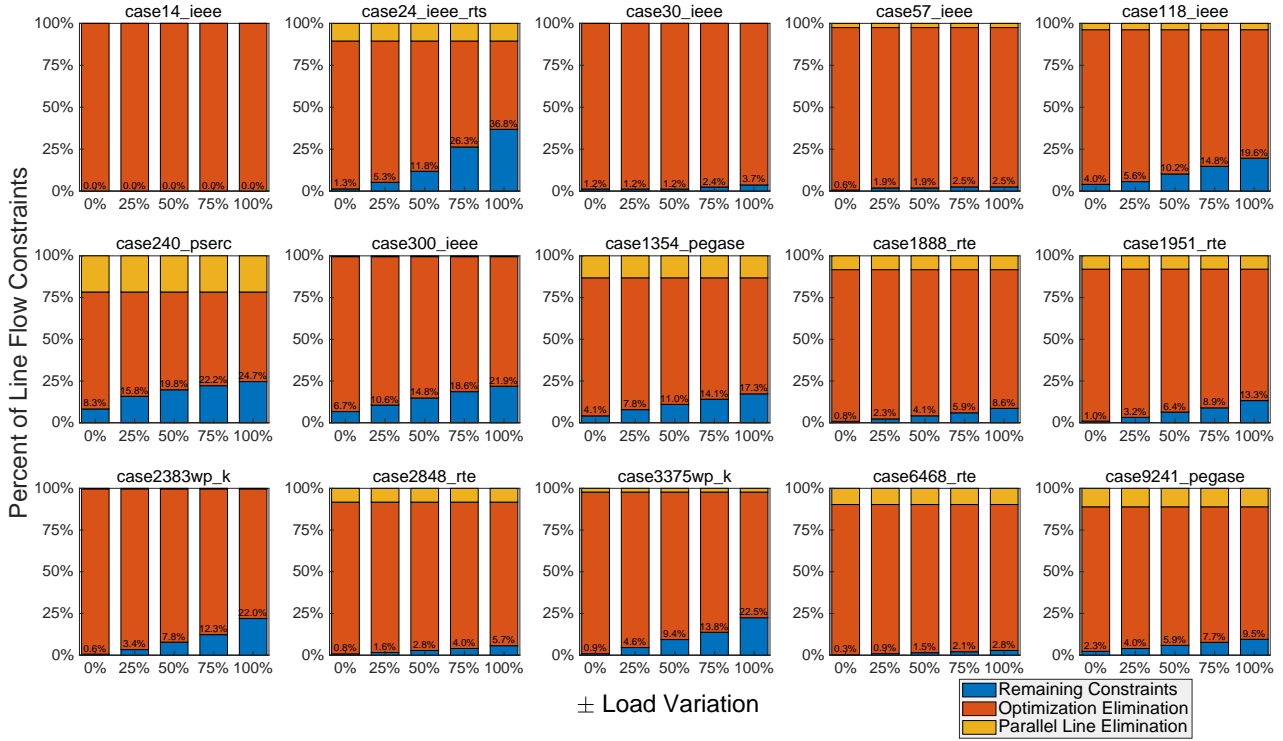}
	\caption{Percent of constraints eliminated by the screening methods.%
	%The values listed below each plot correspond to the number of flow constraints in the original problem.%
	}
	\label{fig:constraint_elim_results}
\end{figure*}

The results for several representative and commonly used test cases are shown in Fig.~\ref{fig:constraint_elim_results}. Each subplot corresponds to one test case. For each test case, we show results for load ranges from $0$\% to $\pm 100$\% of the nominal values, corresponding to the five uncertainty sets described by~\eqref{eq:uncset}. The results of the five load ranges are shown as five bars in each subplot, where each bar represent the percentages of redundant and non-redundant constraints. The yellow portions indicate the percentage of constraints identified as redundant by the parallel line screening, the red portions correspond to the percentage of constraints identified as redundant by the optimization-based constraint screening, and the blue portions denote the percentage of constraints remaining after applying both steps of the screening procedure. For reference, Table~\ref{t:case_summary} summarizes the number of flow constraints in the original (non-screened) problem.

For all test case results shown in Fig.~\ref{fig:constraint_elim_results}, we observe that the fraction of constraints remaining in the problem after the screening is small (i.e., the blue bar represents a small fraction of the total number of constraints). This empirically confirms our claim that the screening methods are effective at identifying redundant constraints for wide ranges of load variation. While the fraction of remaining constraints increases with larger ranges of variation, the screening methods still reduced the number of constraints by between $75.3$\% and $97.2$\% for $\pm 100$\% variation in the load demands. Since the test cases in our simulations correspond to a variety of realistic system models, we expect these results to be representative of many practical power systems.

\begin{table}[t]
\caption{Test Case Summary}
    \centering
\begin{tabular}{|l|c|c|c|}
\hline
     \multicolumn{1}{|c|}{\textbf{Case}} & \textbf{Num.} & \textbf{Num. Flow} & \textbf{UC Solver Time }\\
     \multicolumn{1}{|c|}{\textbf{Name}} & \textbf{Gen.} & \textbf{Constraints} & \textbf{w/o Screening [s]}\\
     \hline\hline
     case14\_ieee & \hphantom{000}5 & \hphantom{000}40 & \hphantom{00}0.006\\
     case24\_ieee\_rts & \hphantom{00}33 & \hphantom{000}76 & \hphantom{00}0.060\\
     case30\_ieee & \hphantom{000}6 & \hphantom{000}82 & \hphantom{00}0.005\\
     case57\_ieee & \hphantom{000}7 & \hphantom{00}160 & \hphantom{00}0.005\\
     case118\_ieee & \hphantom{00}54 & \hphantom{00}372 & \hphantom{00}0.014\\
     case240\_pserc & \hphantom{0}143 & \hphantom{00}896 & \hphantom{00}0.768\\
     case300\_ieee & \hphantom{00}69 & \hphantom{00}822 & \hphantom{00}0.084\\
     case1354\_pegase & \hphantom{0}260 & \hphantom{0}3982 & \hphantom{00}2.753\\
     case1888\_rte  & \hphantom{0}290 & \hphantom{0}5062 & \hphantom{00}1.474\\
     case1951\_rte  & \hphantom{0}366 & \hphantom{0}5192 & \hphantom{00}1.964\\
     case2383wp\_k  & \hphantom{0}327 & \hphantom{0}5792 & \hphantom{00}2.813\\
     case2848\_rte  & \hphantom{0}511 & \hphantom{0}7552 & \hphantom{00}2.413\\
     case3375wp\_k  & \hphantom{0}479 & \hphantom{0}8322 & 30.948\\
     case6468\_rte  & \hphantom{0}399 & 18000 & \hphantom{00}3.710\\
     case9241\_pegase  & 1445 & 32098 & 295.256\\\hline
\end{tabular}
    \label{t:case_summary}
\end{table}

\subsection{Computational Improvements for UC Problems}
\label{subsec:results_uc}
% 2) Computational gains for solving the UC problem

\begin{figure*}[ht]
	\centering
	\includegraphics[width=1.08\textwidth,trim={3.25cm 0cm 1cm 0},clip]{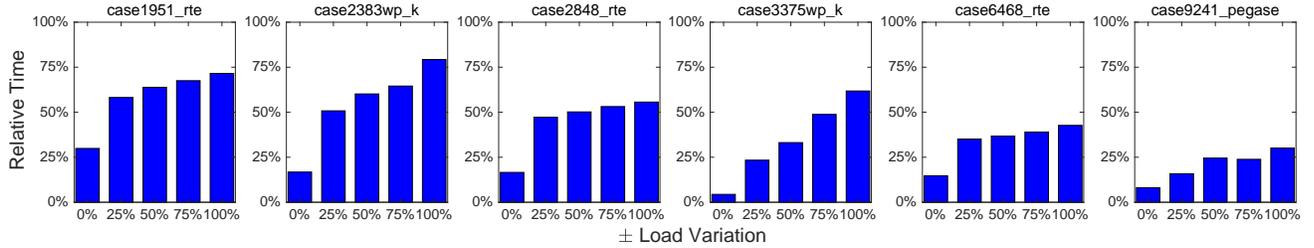}
	\caption{Normalized solver times for the unit commitment problems after constraint screening relative to the solver times prior to constraint screening.%
	%The values listed below each plot correspond to the solver times for the non-screened cases, averaged over all variation amounts.%
	}
	\vspace{-0.2in}
	\label{fig:UC_results}
\end{figure*}

In order to showcase one example of the advantages provided by constraint screening for large variations of load, we consider the single-period UC problem~\eqref{eq:ptdf_uc}. For each test case and each range of load variation, we solve $100$ instances of UC problems~\eqref{eq:ptdf_uc} with load demands uniformly sampled from the corresponding variation range and calculate the average solution time. Fig.~\ref{fig:UC_results} shows representative results for several larger test cases. The blue bars denote the average solution times for the screened problems normalized by the solution times for the original (non-screened) problems. Average solution times in seconds for the original problems are presented in Table~\ref{t:case_summary}. 

Eliminating redundant constraints yields significant computational improvements with solver times between $20.7$\% and $95.6$\% faster than the original problems for the systems with more than 1000 nodes. As expected, the solver times increase with larger ranges of variation as the screening removes fewer redundant constraints. Nevertheless, the screening method still results in substantial computational advantages even for $\pm 100$\% load variation with solver times reductions between $20.7$\% and $69.9$\% over the original problems. While these results do not include the times required for the constraint screening itself, we reiterate that the wide range of variation considered here implies that the screening computations can be conducted offline while remaining applicable to many problems encountered online. Consideration of $\pm 100$\% load variation would make the screening results applicable for typical variations in load over the course of a year~\cite{RTS96}.

We note that while our experiments only considered the single-period DC UC problem~\eqref{eq:ptdf_uc}, the constraints identified as redundant by the screening method remain redundant in the more general multi-period UC problems. %with more general generation constraints.% that are solved in practical applications. 
We expect that constraint screening methods will provide even larger relative improvements in these more complicated problems.

%{\color{red} Finally, we also emphasize that the ability to screen out many constraints despite large ranges of variation suggests the significant potential of constraint screening to address the related applications discussed in Section~\ref{sec:applications}. }

% Improvement expected to get better with more complicated problem descriptions.

\section{Conclusion}
\label{sec:conclusion}

A common observation in power system optimization is that only a limited number of transmission line constraints are ever binding. This paper formalizes this observation using a constraint screening method that rigorously identifies redundant constraints. This method begins by quickly identifying redundancies among parallel lines and then uses an optimization-based method that computes the most extreme achievable values for certain constrained quantities. If the extreme achievable values are less than the specified bounds, the associated constraint is redundant and can be eliminated from the problem. Applying this screening method to a diverse set of large-scale problems reveals that a significant fraction of the flow constraints are redundant even when considering large ranges of variation ($\pm 100$\% of the nominal load demands). An immediate implication of this result is that constraint screening methods can provide computational improvements for a wide range of problems. Importantly, this result also suggests the potential suitability for constraint screening methods in other applications discussed in this paper, 
% (e.g., development of reduced models, stochastic optimization, and ensuring constraint satisfaction in settings with limited observability and control), 
which are subjects of our ongoing work. We are also evaluating the capabilities of constraint screening methods for problems that use the AC power flow model.

\bibliographystyle{IEEEtran}
\bibliography{refs}

\end{document}